\documentclass[11pt, a4paper]{article}
\usepackage[T1]{fontenc}
\usepackage[utf8]{inputenc}
\usepackage{authblk}
\usepackage{amsmath,amsfonts,amssymb,amscd,amsthm,xspace}
\usepackage{algorithm,algpseudocode}
\algrenewcommand\algorithmicrequire{\textbf{Input:}}
\algrenewcommand\algorithmicensure{\textbf{Output:}}
\usepackage{graphicx}
\usepackage{hyperref}
\usepackage{vmargin}
\setmarginsrb  { 0.8in}  
               { 0.6in}  
               { 0.8in}  
               { 0.8in}  
               {  20pt}  
               {0.25in}  
               {   9pt}  
               { 0.3in}  

\title{Performance Analysis of Effective Methods for Solving Band Matrix SLAEs
after Parabolic Nonlinear PDEs}
\author[1]{Milena Veneva \thanks{milena.p.veneva@gmail.com}}
\author[1]{Alexander Ayriyan \thanks{ayriyan@jinr.ru}}
\affil[1]{Joint Institute for Nuclear Research, Laboratory of Information Technologies, Joliot-Curie 6, 141980 Dubna, Moscow region, Russia}

\date{}
\begin{document}
\maketitle
\abstract{%
This paper presents an experimental performance study of implementations of
three different types of algorithms for solving band matrix systems of
linear algebraic equations~(SLAEs) after parabolic nonlinear partial
differential equations -- direct, symbolic, and iterative, the former two of
which were introduced in Veneva and Ayriyan~(arXiv:1710.00428v2). An iterative algorithm is
presented -- the strongly implicit procedure~(SIP), also known as the Stone
method. This method uses the incomplete LU~(ILU(0)) decomposition. An
application of the Hotelling-Bodewig iterative algorithm is suggested as a
replacement of the standard forward-backward substitutions. The upsides and the
downsides of the SIP method are discussed. The complexity of
all the investigated methods is presented. Performance analysis of the
implementations is done using the high-performance computing~(HPC)
clusters ``HybriLIT'' and ``Avitohol''. To that purpose, the experimental setup
and the results from the conducted computations on the individual computer
systems are presented and discussed.
}
\maketitle
%

\section{Introduction}
\label{intro}
Systems of linear algebraic equations~(SLAEs) with pentadiagonal~(PD) and
tridiagonal~(TD) coefficient matrices arise after discretization of partial
differential equations~(PDEs), using finite difference methods~(FDM) or finite
element methods~(FEM). Methods for numerical solving of SLAEs with such
matrices which take into account the band structure of the matrices are
needed. The methods known in the literature usually require the matrix to
possess special characteristics so as the method to be stable, e.g. diagonally
dominance, positive definiteness, etc. which are not always feasible.

In~\cite{ref_name_31}, a finite difference scheme with first-order
approximation of a parabolic PDE was built that leads to a TD SLAE with a
diagonally dominant coefficient matrix. The system was solved using the Thomas method~(see~\cite{ref_name_36}). However, a difference scheme with second-order approximation~\cite{ref_name_30} leads to a matrix which does not have any of
the above-mentioned special characteristics. The numerical algorithms for
solving multidimensional governing equation, using FDM~(e.g. alternating
direction implicit~(ADI) algorithms~(see~\cite{ref_name_43}, \cite{ref_name_35})),
ask for a repeated SLAE solution. This explains the importance of the existence
of effective methods for the SLAE solution stage.

Two different approaches for solving SLAEs with pentadiagonal and tridiagonal
coefficient matrices were explored by us in~\cite{ref_name_30} -- diagonal
dominantization and symbolic algorithms. These approaches led to five
algorithms -- numerical algorithms based on LU decomposition~(for
PD~(see~\cite{ref_name_32}) and TD matrices -- \textbf{NPDM} and \textbf{NTDM}),
modified numerical algorithm for solving SLAEs with a PD
matrix~(where the sparsity of the first and the fifth diagonals was taken into account -- \textbf{MNPDM}), and symbolic algorithms~(for PD~(see~\cite{ref_name_32})
and TD~(see~\cite{ref_name_44}) matrices -- \textbf{SPDM} and \textbf{STDM}). The
numerical experiments with the five methods in our previous paper were
conducted on a PC~(OS: Fedora 25; Processor: Intel Core i7-6700 (3.40 GHz)),
using compiler GCC 6.3.1 and optimization -O0. While the direct numerical
methods have requirements to the coefficient matrix, the direct symbolic
ones only require nonsingularity. Here, we are going to suggest an iterative
numerical method which is also not restrictive on the coefficient matrix.

It is a well-known fact that solving problems of the computational linear
algebra with sparse matrices is crucial for the effectiveness of most of the
programs for computer modelling of processes which are described with the help
of differential equations, especially when solving complex multidimensional
problems. However, this is exactly how most of the computational science
problems look like and hence usually they cannot be modelled on ordinary PCs for
a reasonable amount of time. This enforces the usage of supercomputers and
clusters for solving such big problems. For example, a numerical solving of a
parabolic PDE needs to solve independently (or in parallel) $N$ SLAEs $d$ times
at each time-step, where $N$ is the discretization number, i.e. the matrix
dimension, and $d$ is the dimension of the PDE. Thus, it is also important to
have an efficient method for serial solving of one band SLAE. Therefore, the aim
and the main contribution of this paper is to investigate the performance
characteristics of the considered serial methods for band SLAE being executed on
modern computer clusters.

The layout of the paper is as follows: in the next section, we introduce the
outline of the SIP algorithm. Afterwards, we introduce the experimental setup
including the description of the computers used in our experiments and analyze
the obtained results.


\section{Iterative Approach}
\label{sec-1}

An iterative procedure for solving SLAEs with a pentadiagonal coefficient matrix
is considered, namely the strongly implicit procedure~(SIP)~(see~\cite{ref_name_39}),
also known as the Stone method. It is an algorithm for solving sparse SLAEs.
The method uses the incomplete LU~(ILU(0)) decomposition~(see~\cite{ref_name_40})
which is an approximation of the exact LU decomposition in the case when a
sparse matrix is considered. The idea of ILU(0) is that the zero elements of
$L$ and $U$ are chosen to be on the same places as of the initial matrix $A$.
In the case of a pentadiagonal coefficient matrix $A$, $LU$ is going to be also
pentadiagonal, $L$ and $U$ are going to have non-zero elements only on three of
their diagonals (main diagonal and two subdiagonals for $L$; main diagonal and
two superdiagonals for $U$). The Stone method for solving a SLAE of the form
$Ax=b$ can be seen in Algorithm~\ref{alg:1}. There, $LU$ is found using the
ILU(0) algorithm suggested in~\cite{ref_name_40}; $L$ and $U$ are extracted
using a modification of the Doolittle method~(see~\cite{ref_name_43}), namely
instead of referencing the matrix $A$, we reference the already found $LU$
matrix. This way the product of $L$ and $U$ is exactly $LU$.
\begin{algorithm}
\caption{The Stone method for solving a SLAE $Ax = b$}
\label{alg:1}
\begin{algorithmic}[1]
\Require{$A, b, L, U, errorMargin$}
\Ensure{$k, \overrightarrow{x}^{(k)}$}
\State{$k=0$} \Comment{number of iterations}
\State{$\overrightarrow{x}^{(k)} = \overrightarrow{0}$} \Comment{set an initial guess vector}
\State{$\overrightarrow{newRHS}^{(k)} = A\,\overrightarrow{x}^{(k)}$} \Comment{new right-hand side~(RHS)}
\State{$\overrightarrow{residual}^{(k)} = \overrightarrow{b} - \overrightarrow{newRHS}^{(k)}$}
\State{$K = L\,U - A$}
\While{$\|\overrightarrow{residual}^{(k)}\|_{\infty}\geq\textrm{errorMargin}$}
\State{$\overrightarrow{newRHS}^{(k)} = K\,\overrightarrow{x}^{(k)} + \overrightarrow{b}$}
\State{solve $L\,\overrightarrow{y}^{(k)} = \overrightarrow{newRHS}^{(k)}$}\label{row:1}
\State{solve $U\,\overrightarrow{x}^{(k+1)} = \overrightarrow{y}^{(k)}$}\label{row:2}
\State{$\overrightarrow{residual}^{(k+1)} = \overrightarrow{b} - A\,\overrightarrow{x}^{(k+1)}$}
\State{$k++$}
\EndWhile
\end{algorithmic}
\end{algorithm}
\noindent Every iteration step of the Stone method consists of two
matrix-vector multiplications with a pentadiagonal matrix, one forward and one
backward substitutions with the two triangular matrices of the ILU(0), and two
vector additions, i.e. the complexity of the algorithm on every iteration is
$31\,N-36 = O(N)$, where $N$ is the number of rows of the initial matrix.

\noindent{\textbf{Remark:} Instead of using forward and backward substitutions
on rows \ref{row:1}-\ref{row:2} of Algorithm~\ref{alg:1}, one can try to
find the inverse matrices of $L$ and $U$, using a numerical procedure,
e.g.\,the Hotelling-Bodewig iterative algorithm~(see~\cite{ref_name_41}). A diagonal
matrix can be used as an initial guess for the inverse matrix, as it is
suggested in~\cite{ref_name_42}. Since a matrix implementation is going to be
very demanding in regards to memory, conduction of computational experiments
for a matrix with more than $7\times10^3$ rows is going to be impossible. For
that reason, the algorithms could be redesigned, taking into account the band
structure of the data, and so an array implementation could be made. (For the Hotelling-Bodewig iterative algorithm and numerical results from that approach,
see Appendix~\ref{appendix:1}.)


\section{Numerical Experiments}
\label{sec-3}

Computations were held on the basis of the heterogeneous computing cluster
``HybriLIT'' at the Laboratory of Information Technologies of the Joint
Institute for Nuclear Research in town of science Dubna, Russia and on the
cluster computer system ``Avitohol'' at the Advanced Computing and Data
Centre of the Institute of Information and Communication Technologies of the
Bulgarian Academy of Sciences in Sofia, Bulgaria.

\subsection{Experimental Setup}
The direct and iterative numerical algorithms are implemented using \texttt{C++},
while the symbolic algorithms are implemented using the \texttt{GiNaC}
library~(version 1.7.2) (see~\cite{ref_name_33}) of \texttt{C++}.

The heterogeneous computing cluster ``HybriLIT'' consists of 13 computational
nodes which include two Intel Xeon E5-2695v2 processors (12-core) or
two Intel Xeon E5-2695v3 processors (14-core). For more information, visit\\
\url{http://hybrilit.jinr.ru/en}. It must be mentioned that for the sake
of the performance analysis and the comparison between the computational times
only nodes with Intel Xeon E5-2695v2 processors were used.

The supercomputer system ``Avitohol'' is built with HP Cluster Platform SL250S
GEN8. It has two Intel 8-core Intel Xeon E5-2650 v2 8C processors each of which
runs at 2.6 GHz. For more information, visit \url{http://www.hpc.acad.bg/}.
``Avitohol'' has been part of the TOP500 list (\url{https://www.top500.org})
twice -- ranking 332nd in June 2015 and 388th in November 2015.

Tables~\ref{tab:01} and \ref{tab:02} summarize the basic information about
hardware, compilers and libraries used on the two computer systems. The
reason why different compilers were used for the numerical and the
symbolic methods, respectively, is that the \texttt{GiNaC} library does not
maintain work with the \texttt{Intel} compilers. However, for the numerical methods
the \texttt{Intel} compilers gave us better results than the \texttt{GCC} ones.
\renewcommand{\arraystretch}{1.2}
\begin{table}[!htb]
\centering
\begin{tabular}{ccc}
\hline\hline
Computer system & Processor & Number of processors per node \\\hline\hline
``HybriLIT''	&  Intel Xeon E5-2695v2 & 2 \\
                &  Intel Xeon E5-2695v3 & 2 \\\hline\hline
``Avitohol''    &  Intel Xeon E5-2650v2 & 2 \\\hline\hline       
\end{tabular}
\caption{Information about the available hardware on the two computer systems}
\label{tab:01}
\end{table}
\renewcommand{\arraystretch}{1.2}
\begin{table}[!htb]
\centering
\begin{tabular}{r|c|c}
\hline\hline
Computer                    & ``HybriLIT''          & ``Avitohol''   \\
system                      &                       &            \\\hline\hline
Compiler for the direct     & Intel 2017.2.050 ICPC & Intel 2016.2.181 ICPC \\
and iterative procedures    &                       &            \\\hline\hline
Compiler for the            & GCC 4.9.3             & GCC 6.2.0  \\
symbolic procedures         &                       &            \\\hline\hline
Needed libraries for the    & GiNaC (1.7.2)         & GiNaC (1.7.2) \\
symbolic procedures         & CLN (1.3.4)           & CLN (1.3.4) \\\hline\hline
Optimization for the direct & -O2                   & -O2         \\
and iterative procedures    &                       &             \\\hline\hline
Optimization for the        & -O0                   & -O0         \\
symbolic procedures         &                       &             \\\hline\hline
\end{tabular}
\caption{Information about the used software on the two computer systems}
\label{tab:02}
\end{table}

\subsection{Experimental Results}

During our experiments wall-clock times were collected using the
member function \texttt{now()} of the class
\texttt{std::chrono::high\_resolution\_clock} which represents
the clock with the smallest tick period provided by the implementation; it
requires at least standard \texttt{c++11} (needs the argument -std=c++11 when
compiling). We report the average time from multiple runs. Since the largest
supported precision in the \texttt{GiNaC} library is \textbf{double}, during
all the experiments double data type is used. The achieved accuracy during all
the numerical experiments is summarized, using infinity norm. The notation is
as follows: \textbf{NPDM} stands for numerical~PD method, \textbf{MNPDM} --
modified numerical~PD method, \textbf{SPDM} -- symbolic~PD method,
\textbf{NTDM} -- numerical~TD
method, \textbf{STDM} -- symbolic~TD method. The error tolerance used in
the iterative method is $10^{-12}$. Both the methods comprised in the iterative
procedure (ILU(0) and SIP) are implemented using an array representation of the
matrices instead of a matrix one.

\noindent\textbf{Remark 1:} So as the nonsingularity of the matrices to be
checked, a fast symbolic algorithm for calculating the determinant is
implemented, using the method suggested in~\cite{ref_name_37}. The
complexity of the algorithm is $O(N)$.

\noindent\textbf{Remark 2:} The number of needed operations for the Gaussian
elimination used so as~PD matrices to be transformed into~TD ones is $18+16\,K$,
where $K$ is the number of~PD matrix rows with nonzero elements on their second
subdiagonal and on their second superdiagonal. Usually, $K\ll N$.

%

The achieved computational times from solving a SLAE on ``HybriLIT'' are
summarized in Tables~\ref{tab:h1} and \ref{tab:h2}. The number of needed
iterations for the Stone method is 31.
\begin{table}[!htb]
\centering
\begin{tabular}{crrrrr}
\hline\hline
& \multicolumn{5}{c}{Wall-clock time [s]}  \\\hline\hline
$N$      & NPDM 	& MNPDM    & SPDM        & NTDM 	  & STDM\\\hline
$10^{3}$ & 0.0000427& 0.0000420&  0.1098275  & 0.0000273  & 0.0827742\\
$10^{4}$ & 0.0004310& 0.0004270& 17.5189275  & 0.0002693  & 7.9570979\\
$10^{5}$ & 0.0041760& 0.0040850&5991.4962896 & 0.0026823  & 2857.5483843\\
$10^{8}$ & 2.7946627& 2.6662850& \multicolumn{1}{c}{--} & 2.0525187 & \multicolumn{1}{c}{--} \\
\hline\hline
$\max\limits_{\textrm{\scriptsize{N}}}\|y - \bar{y}\|_{\infty}$ &
\multicolumn{1}{c}{$2.22\times 10^{-16}$} &
\multicolumn{1}{c}{$2.22\times 10^{-16}$} &
\multicolumn{1}{c}{0} &
\multicolumn{1}{c}{$2.22\times 10^{-16}$} &
\multicolumn{1}{c}{0} \\
\hline\hline
\end{tabular}
\caption{Results from solving a SLAE on the cluster ``HybriLIT'' applying direct methods.}
\label{tab:h1}
\end{table}

\begin{table}[!htb]
\centering
\begin{tabular}{crr}
\hline\hline
& \multicolumn{2}{c}{Wall-clock time [s]}                             \\\hline\hline
$N$            & \multicolumn{1}{c}{ILU(0)} & \multicolumn{1}{c}{SIP} \\\hline
$10^2$         &   0.0004090                &  0.0000847              \\
$10^3$         &   0.3493933                &  0.0007300              \\
$10^4$         & 325.6368140                &  0.0084683              \\\hline\hline
$\max\limits_{\textrm{\scriptsize{N}}}\|y - \bar{y}\|_{\infty}$ &
\multicolumn{1}{c}{--} &
\multicolumn{1}{c}{$1.42\times 10^{-13}$}\\
\hline\hline
\end{tabular}
\caption{Results from solving a SLAE using SIP on the cluster ``HybriLIT'' applying an iterative method.}
\label{tab:h2}
\end{table}
%
%

Tables ~\ref{tab:a1} and \ref{tab:a2} sum up the computational times from
solving a SLAE on ``Avitohol''. The number of needed iterations for the Stone
method is 31.
\begin{table}[!htb]
\centering
\begin{tabular}{crrrrr}
\hline\hline
& \multicolumn{5}{c}{Wall-clock time [s]}  \\\hline\hline
$N$      & NPDM 	& MNPDM      & SPDM 	 & NTDM 	 & STDM\\\hline
$10^{3}$ & 0.0000420 & 0.0000400 & 0.1089801 & 0.0000290 & 0.0518055 \\
$10^{4}$ & 0.0004234 & 0.0004110 & 15.2383414& 0.0002610 & 5.2806483 \\
$10^{5}$ & 0.0040710 & 0.0039387 &2009.6004854& 0.0027417 & 711.9402796 \\
$10^{8}$ & 2.8660797 & 2.7304760 & \multicolumn{1}{c}{--}& 2.1347700   & \multicolumn{1}{c}{--}\\
\hline\hline
$\max\limits_{\textrm{\scriptsize{N}}}\|y - \bar{y}\|_{\infty}$ &
\multicolumn{1}{c}{$2.22\times 10^{-16}$} &
\multicolumn{1}{c}{$2.22\times 10^{-16}$} &
\multicolumn{1}{c}{0} &
\multicolumn{1}{c}{$2.22\times 10^{-16}$} &
\multicolumn{1}{c}{0} \\
\hline\hline
\end{tabular}
\caption{Results from solving a SLAE on the supercomputer system ``Avitohol'' applying direct methods.}
\label{tab:a1}
\end{table}
\vspace{-0em}
\begin{table}[!htb]
\centering
\begin{tabular}{crr}
\hline\hline
& \multicolumn{2}{c}{Wall-clock time [s]}                             \\\hline\hline
$N$            & \multicolumn{1}{c}{ILU(0)} & \multicolumn{1}{c}{SIP} \\\hline
$10^2$         &   0.0011817                & 0.0001210               \\
$10^3$         &   0.5288667                & 0.0008603               \\
$10^4$         &   516.6088950              & 0.0085333               \\\hline\hline
$\max\limits_{\textrm{\scriptsize{N}}}\|y - \bar{y}\|_{\infty}$ &
\multicolumn{1}{c}{--} &
\multicolumn{1}{c}{$1.42\times 10^{-13}$}\\
\hline\hline
\end{tabular}
\caption{Results from solving a SLAE using SIP on the cluster ``Avitohol'' applying an iterative method.}
\label{tab:a2}
\end{table}


\section{Discussion and Conclusions}
\label{sec-4}

Three different approaches for solving a SLAE are compared -- direct numerical,
direct symbolic, and iterative. The complexity of all the suggested numerical
algorithms is $O(N)$ (see Table~\ref{tab:301}). Since it is unknown what stands
behind the symbolic library, evaluating the complexity of the symbolic algorithms
is a very complicated task.
\begin{table}[!htbp]
\centering
\begin{tabular}{c|c|c|c|c|c|c}
\hline\hline
Method:     & NPDM 	     & MNPDM       & SPDM & NTDM   & STDM & SIP      \\\hline
Complexity: & $19N - 29$ & $13N+7K-14$ & --   & $9N+2$ & --   & $31N-36$ \\\hline\hline
\end{tabular}
\caption{Complexity of the investigated methods.}
\label{tab:301}
\end{table}

Both the achieved computational times and accuracy for the \textbf{NPDM} and
\textbf{SPDM} methods on both the clusters were much better than the ones
outlined in~\cite{ref_name_32}.

All the experiments with the direct methods gave an accuracy of an order of
magnitude of $10^{-16}$, while the iterative method gave an accuracy of an
order of magnitude of $10^{-13}$.
\begin{figure}[!htb]
  \centering
  \includegraphics[width=\textwidth]{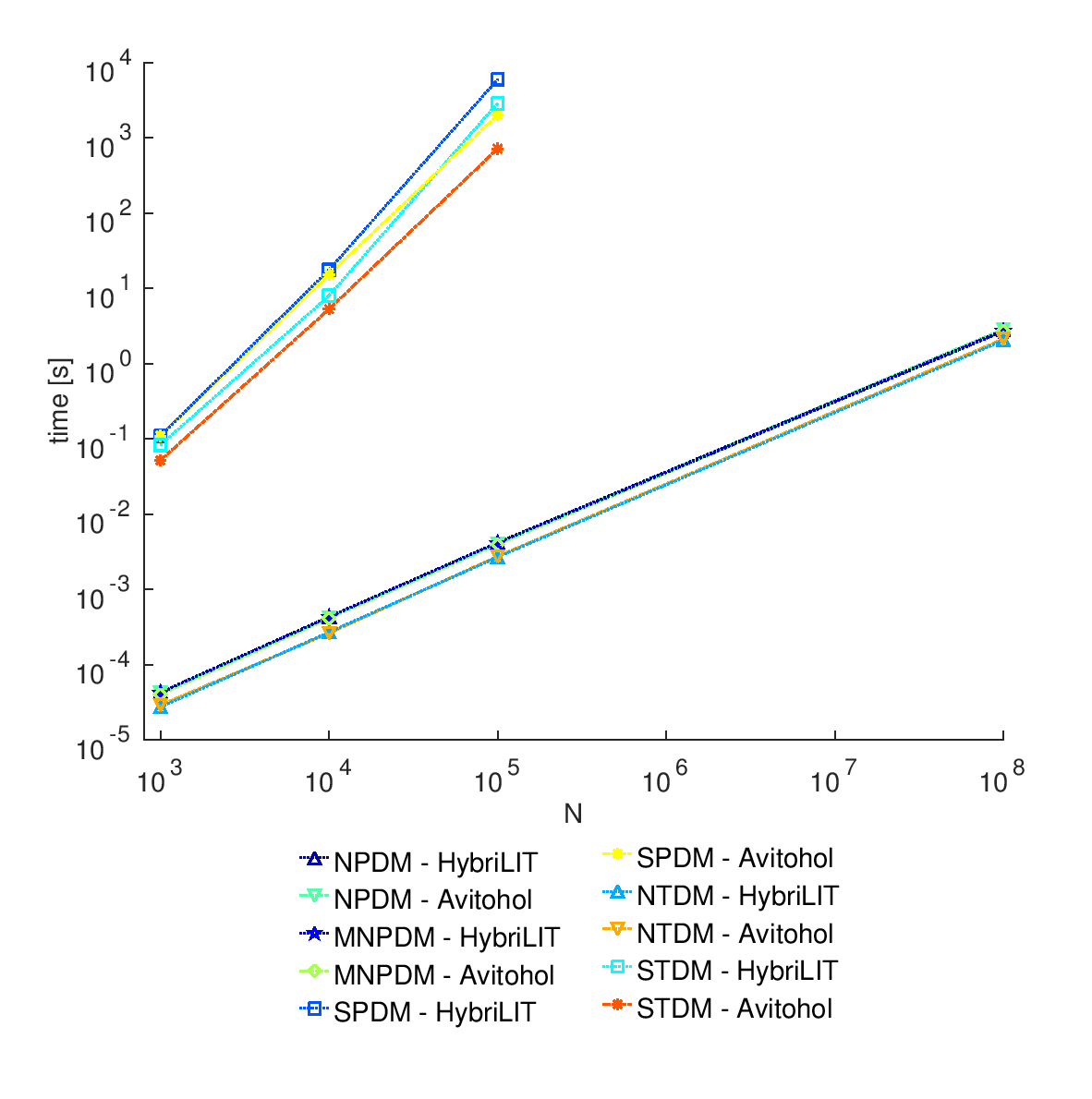}
  \caption{Execution time comparison.}
  \label{fig:1}
\end{figure}
Expectedly, the modified version of the numerical method for solving a SLAE with
a~PD matrix \textbf{MNPDM} gave better computational time than the general
algorithm \textbf{NPDM} in the case of a sparse~PD coefficient matrix, since
the former method has a lower complexity~(usually $K\ll N$, where $K$ is the
number of~PD matrix rows with nonzero elements on their second subdiagonal
and on their second superdiagonal). The fastest numerical
algorithm was found to come from the Thomas method. Finally, an iterative
algorithm was built -- the Stone method. For the needs of the method,
additionally, ILU(0) was implemented. An upside of this iterative
procedure is that it requires the initial matrix to be nonsingular only.
However, this method is not suitable for matrices for which
$N > 1\times10^5$, since the ILU(0) decomposition of a matrix is
computationally demanding on time and memory. Here, likewise the
symbolic algorithms, in the case of a piecewise linear parabolic
partial differential equation, they do not add nonlinearity to the
right-hand side of the system and hence, there is no need of
iterations for the time step to be executed~(see~\cite{ref_name_30}). Similarly
to the symbolic methods, SIP is not comparable with the numerical algorithms with
respect to the required time in the case of a numerical solving of the heat
equation when one needs to solve the SLAE many times. Lastly, the obtained
accuracy is worse in comparison with any of the other methods.
A comparison between the execution times for the direct numerical
methods (see Figure~\ref{fig:1}) showed only a negligible
difference between the two computer systems.
However, this was not the case when it comes to the symbolic methods where
``Avitohol'' performed much better than ``HybriLIT''. On the other hand,
``HybriLIT'' behaved better than ``Avitohol'' with respect to the ILU(0)
procedure (see Figure~\ref{fig:2}). Only a minimal discrepancy
in times was observed for the SIP algorithm.
\begin{figure}[!htb]
  \centering
  \includegraphics[width=\textwidth]{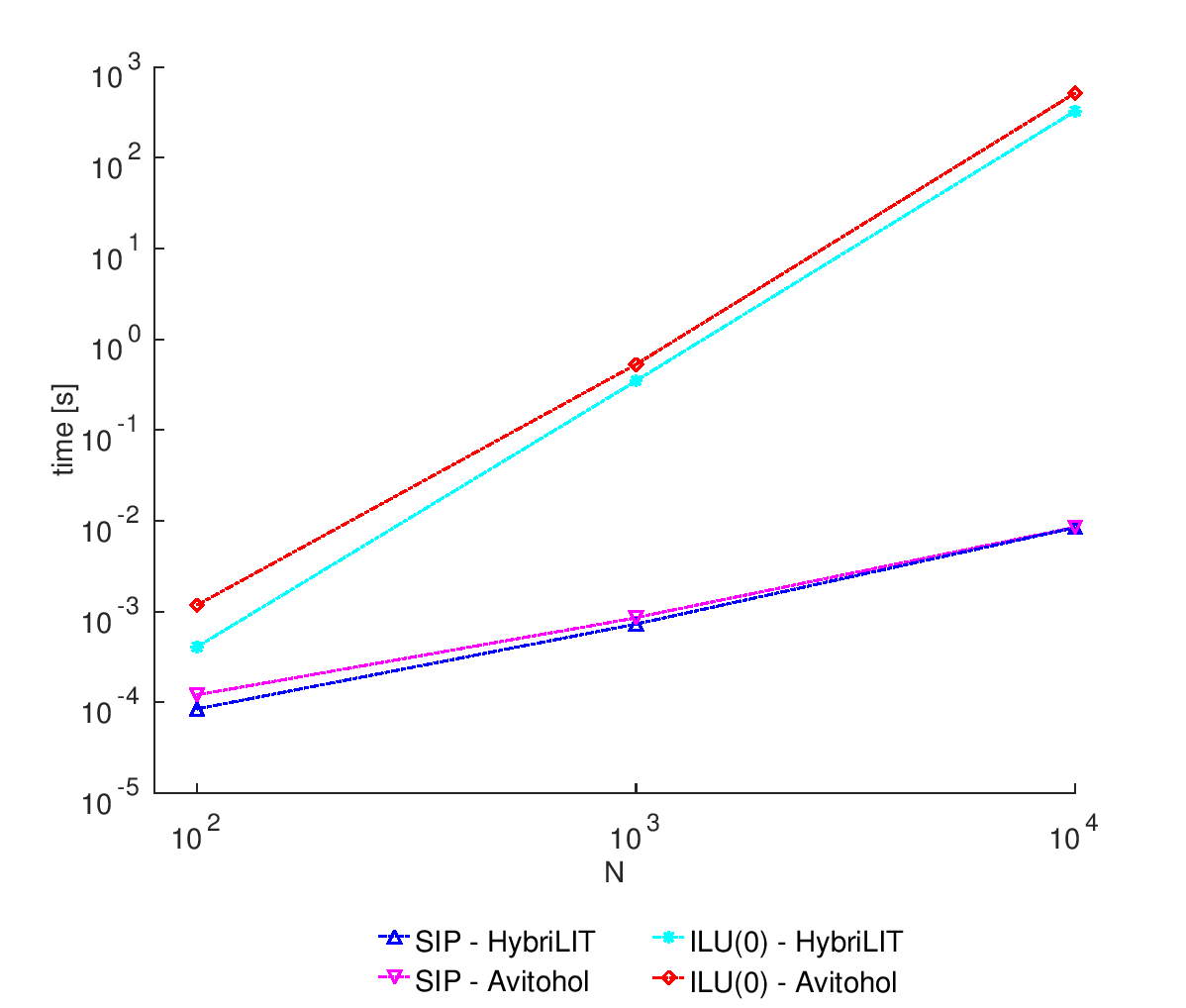}
  \caption{Execution time comparison.}
  \label{fig:2}
\end{figure}


\section*{Acknowledgements}
The authors want to express their gratitude to the Summer Student Program at
JINR, Dr.\,J\'{a}n Bu\v{s}a Jr. (JINR), Dr.\,Andrey Lebedev (GSI/JINR), Assoc.\,
Prof.\,Ivan Georgiev (IICT \& IMI, BAS), the ``HybriLIT'' team at LIT, JINR, and
the ``Avitohol'' team at the Advanced Computing and Data Centre of IICT, BAS.
Computer time grants from LIT, JINR and the Advanced Computing and Data Centre at
IICT, BAS are kindly acknowledged.
A. Ayriyan thanks the JINR grant No. 17-602-01.


\section*{Appendix}
\label{appendix:1}
%
The Hotelling-Bodewig iterative algorithm has the form as follows:
\begin{equation}
\label{eq:3_9}
A^{-1}_{n+1} = A^{-1}_{n}\,(2\,I - A\,A^{-1}_{n}), \quad n=0,1,\ldots,
\end{equation}
\noindent where $I$ is the identity matrix, $A$ is the matrix whose inverse we
are looking for. $A^{-1}_{0}$ is taken to be of a diagonal form.

The obtained computational times for the ILU(0) method, the
Hotelling-Bodewig iterative algorithm and the Stone method, using the
heterogeneous cluster ``HybriLIT'' and the supercomputer system ``Avitohol'',
are summarized in Tables~\ref{tab:app1}, \ref{tab:app2}, and \ref{tab:app3}.
\begin{table}[!htb]
\centering
\begin{tabular}{crrr||rrr}
\hline\hline
& \multicolumn{6}{c}{Wall-clock time [s]} \\\hline\hline
& \multicolumn{3}{c||}{matrix implementation} & \multicolumn{3}{c}{array implementation} \\\hline\hline
$N$ & \multicolumn{1}{c}{ILU(0)} & \multicolumn{1}{c}{$L^{-1}$} & \multicolumn{1}{c||}{$U^{-1}$} & \multicolumn{1}{c}{ILU(0)}& \multicolumn{1}{c}{$L^{-1}$} & \multicolumn{1}{c}{$U^{-1}$}\\\hline
$10^2$         &  0.0007027 &  0.0071077 & 0.0096933
               &  0.0004687 &  0.0026583 & 0.0043893 \\
$10^3$         &  1.5635590 & 82.9383600 &49.5268260
               &  0.3368320 &  3.3851580 & 5.3989410 \\
$2\times 10^3$ & 21.6416160 &289.3253220 &300.0902740
               &  2.5914510 & 27.4874390 & 41.7962950\\
$5\times 10^3$ &547.7717120 &4835.9211180&6800.0948670
               & 39.5945850 &1153.8804500&1606.9331050\\
$7\times 10^3$ &1178.6338560&18966.0135900&24345.2476050
               &108.1988910 &3395.9828320&7116.1639450 \\
$10^4$         & \multicolumn{1}{c}{--}  & \multicolumn{1}{c}{--} & \multicolumn{1}{c||}{--}
               &314.7906570 &10384.6694270&14561.3854660 \\
\hline\hline
\end{tabular}
\caption{Results from the ILU(0) method and the numerical method for inverting matrices, using the cluster ``HybriLIT''.}
\label{tab:app1}
\end{table}
\begin{table}[!htb]
\centering
\begin{tabular}{crrr||rrr}
\hline\hline
& \multicolumn{6}{c}{Wall-clock time [s]} \\\hline\hline
& \multicolumn{3}{c||}{matrix implementation} & \multicolumn{3}{c}{array implementation} \\\hline\hline
$N$ & \multicolumn{1}{c}{ILU(0)} & \multicolumn{1}{c}{$L^{-1}$} & \multicolumn{1}{c||}{$U^{-1}$} & \multicolumn{1}{c}{ILU(0)}& \multicolumn{1}{c}{$L^{-1}$} & \multicolumn{1}{c}{$U^{-1}$}\\\hline
$10^2$         &   0.0017620 &    0.0089710 &    0.0117817
               &   0.0013103 &    0.0035383 &    0.0060527 \\
$10^3$         &   1.9317270 &   85.0694670 &   76.6738290
               &   0.5320370 &    4.9676690 &    6.5183100 \\
$2\times 10^3$ &  27.3982830 &  299.0649410 &  370.7769350
               &   4.1901280 &   32.8010570 &   51.0338640  \\
$5\times 10^3$ & 495.6995570 & 5175.7197290 & 6720.2701160
               &  64.8352820 & 1227.5802660 & 1780.3281350 \\
$7\times 10^3$ &1144.9877790 &14829.3973560 &22415.4835190
               & 177.4153890 & 3569.6018970 & 5279.6295710 \\
$10^4$         & \multicolumn{1}{c}{--} & \multicolumn{1}{c}{--} & \multicolumn{1}{c||}{--}
               & 516.4751790 &10441.5862030 &17833.2337200 \\
\hline\hline
\end{tabular}
\caption{Results from the ILU(0) method and the numerical method for inverting matrices, using the cluster ``Avitohol''.}
\label{tab:app2}
\end{table}
\begin{table}[!htb]
\centering
\begin{tabular}{cr|r||r|r}
\hline\hline
& \multicolumn{4}{c}{Wall-clock time [s]}  \\\hline\hline
& \multicolumn{2}{c||}{on ``HybriLIT''} & \multicolumn{2}{c}{on ``Avitohol''} \\\hline\hline
& \multicolumn{1}{c|}{matrix} & \multicolumn{1}{c||}{array} & \multicolumn{1}{c|}{matrix} & \multicolumn{1}{c}{array}\\
& \multicolumn{1}{c|}{implementation} & \multicolumn{1}{c||}{implementation} & \multicolumn{1}{c|}{implementation} & \multicolumn{1}{c}{implementation} \\\hline\hline
$N$ & \multicolumn{1}{c|}{SIP} & \multicolumn{1}{c||}{SIP} & \multicolumn{1}{c|}{SIP} & \multicolumn{1}{c}{SIP} \\\hline
$10^2$         & 0.0005637 & 0.0001827 & 0.0006510 & 0.0002953 \\
$10^3$         & 0.0703420 & 0.0130150 & 0.0866850 & 0.0163560 \\
$2\times 10^3$ & 0.3403310 & 0.0683440 & 0.3492530 & 0.0859700 \\
$5\times 10^3$ & 2.3330770 & 0.5063490 & 3.7949870 & 0.5812540 \\
$7\times 10^3$ & 8.7838330 & 1.1616650 & 6.3790020 & 1.2195610 \\
$10^4$         & \multicolumn{1}{c|}{--} & 2.0574280 & \multicolumn{1}{c|}{--} & 2.9845790 \\\hline\hline
$\max\limits_{\textrm{\scriptsize{N}}}\|y - \bar{y}\|_{\infty}$ &
\multicolumn{1}{c|}{$3.13\times 10^{-14}$} &
\multicolumn{1}{c||}{$3.13\times 10^{-14}$} &
\multicolumn{1}{c|}{$3.13\times 10^{-14}$} &
\multicolumn{1}{c}{$3.13\times 10^{-14}$}\\
\hline\hline
\end{tabular}
\caption{Results from solving a SLAE using SIP on the clusters ``HybriLIT'' and ``Avitohol''.}
\label{tab:app3}
\end{table}

The matrix implementations lead to 5, 7, and 34 iterations, respectively for
finding $L^{-1}$ and $U^{-1}$, applying the Hotelling-Bodewig procedure, and
for the Stone method while the needed iterations when the array
implementations are executed are 5, 6, and 31, respectively. It is expected
that inverting $L$ would require less number of iterations, since it is a
unit triangular matrix. The achieved
accuracy is of an order of magnitude of $10^{-13}$, having used an error
tolerance $10^{-12}$. Comparing the  results for the computational times, one
can see that the array implementation not only decreased the time needed for
the inversion of both the matrices $L$ and $U$ but also it decreases the
number of iterations needed so as the matrix $U$ to be inverted. As one can
see, the time required for the SIP procedure is also improved by the new
implementation approach. One reason being is that the number of iterations is
decreased. Overall, the array implementations decrease the computational times
with one order of magnitude. Finally, this second approach requires less amount
of memory (instead of keeping $N\times N$ matrix, just 5 arrays with length $N$
are stored), which allows experiments with bigger matrices to be conducted.
However, this method (even in its array form) is not suitable for too large
matrices (with number of rows bigger than $1\times10^5$), since the evaluation
of the inverse of a matrix is computationally demanding on both time and memory.
A comparison between the times on the two computer systems showed
that overall ``HybriLIT'' is a bit faster than ``Avitohol''.


\end{document}